\documentclass[10pt]{article}

\usepackage[utf8]{inputenc}
\usepackage[english]{babel}

\usepackage{amsmath, amssymb, bm}

\usepackage{graphicx}
\usepackage{booktabs}
\usepackage{multirow}
\usepackage{makecell}
\usepackage{enumitem}

\usepackage[margin=1.5in]{geometry}
\makeatletter
\renewcommand\footnoterule{%
  \vspace{0.7em}   
  \hrule width 0.4\columnwidth height 0.4pt
  \kern 0.3em      
}
\makeatother
\usepackage[hang]{footmisc}
\usepackage[backend=biber,style=apa]{biblatex}
\usepackage[colorinlistoftodos]{todonotes}

\begin{document}

\title{\textbf{Unique Preference Aggregation in Design and Decision Making}}

\author{
A.R.M. Wolfert\thanks{Corresponding author: A.R.M. (Rogier) Wolfert: \texttt{a.r.m.wolfert@tudelft.nl}}\\[1.2em]  
\small Department of Algorithmics, Faculty of Mathematics and Computer Science\\[-0.3em]  
\small Delft University of Technology, The Netherlands
}

\date{}
\maketitle

\begin{abstract}
\noindent
Preference aggregation is a core operation in multi-objective design optimisation and group decision-making, as it determines the best-fit-for-common-purpose alternative within complex socio-technical contexts. Since preferences are intrinsically linked to choice, they are subjective, inherently contextual, and reflect humans’ will to relatively order alternatives across different criteria. Their aggregation therefore requires a rigorous measurement-theoretic foundation to ensure mathematical validity, interpretability, and uniqueness. Barzilai’s Preference Function Modelling (PFM) theory provides such a foundation by representing preferences as elements of a one-dimensional affine space, in which only differences are meaningful. While PFM establishes the axiomatic conditions for valid preference aggregation, it does not, in itself, provide an operational aggregation procedure for applied design and decision-making.

\noindent
In this paper, a unique PFM-based preference aggregation approach is operationalised for multi-criteria decision-making (MCDM). To this end, a linear preference space is constructed with commensurability across criteria, a stable zero-referenced interval scale, and uniqueness up to affine transformation. It is shown that only linear aggregation—a weighted centroid of \(z\)-normalised preference scores—is admissible. Analytical arguments and illustrative examples demonstrate that this approach yields a single, consistent, and uniquely determined aggregated preference ranking. 

\noindent
Furthermore, it is shown that commonly used aggregation approaches in MCDM—such as weighted arithmetic and geometric means, as well as distance-based optimisation methods—may fail to produce consistent and preference-meaningful rankings, and therefore cannot guarantee decision-valid outcomes. In contrast, the PFM-based aggregation developed in this work provides a principled and implementable framework for robust multi-criteria decision analysis (MCDA) and multi-objective design optimisation (MODO), thereby bridging the gap between measurement theory and applied group design and decision-making in complex problems.

\end{abstract}

\vspace*{1.5em}

\noindent\textbf{Keywords:}
Multi-Criteria Decision Making (MCDM);
Preference Aggregation;
Design and Decision Systems;
Multi-criteria Decision Analysis (MCDA);
Multi-objective Design Optimization (MODO);
Preference Function Modelling (PFM)

\vspace{0.75em}

\noindent\textbf{Article history:}\\
Compiled March~18,~2026

\newpage
\section*{Introduction}

\textbf{Preference} is the \textbf{decisive quantity} in engineering design optimisation and management science for multi-criteria decision-making (MCDM). A preference expresses the relative desirability, \textit{value}, or \textit{utility} of a design alternative or decision option $A_i$ with respect to a criterion $C_j$. Everything of value is relative. Each alternative or option is intuitively evaluated against one’s conscious lived experience --- a relative, subjective, and open-ended human perception arising from all the outer and inner senses. Preference is not a physical property but a subjective construct of the mind. It represents an individual’s freedom of choice — free will — within the set of available options, defining the decision space from which selections are made. Free choice cannot be measured directly like a physical property, because it is not an object of thought but a subjective reality expressed through human willing; what can be assessed is the relative ordering of alternatives, in which preference emerges as a comparative expression of value. The source of value is not the object itself but the relationship between the individual and the world. Preference is inherently contextual — it reflects the free ordering of alternatives within a given situation. It is therefore individual, relational, and situation-dependent. Preference is not an attribute of an alternative but of a comparison. Its existence presupposes a perceived difference between alternatives relative to a reference. Hence preference is fundamentally a relational-relative quantity: without perceived difference there is no preference, and without preference there is no choice. This relational-relative nature is reflected in its affine measurement structure.\\
Consequently, \textbf{preference induces choice}: it is a binary relation over alternatives that determines selection, whereby one alternative is chosen over another: i.e., \(A_1 \succ A_2\). Without differences, no decision — only differences can activate choice. Preference scores or ratings are points whose meaning is inherently relational, as they have no absolute zero. Thus, they are elements of a one-dimensional affine space rather than a ratio or absolute scale (e.g., like time without a fixed origin, a rod without a defined length unit, or a road without kilometre markers — in all cases, only differences and ratios of differences carry meaning). They represent orderings and relative magnitudes of preference differences, but are not physical or absolute measurements \cite{french2023reflections}. Thus, even a seemingly absolute \(0\)--\(100\) scale expressing raw preference ratings for a single criterion remains a \textbf{one-dimensional affine preference scale} (\cite{Barzilai2022}), whose numerical values are defined only up to an affine transformation.

Only ratios of differences between preference values are meaningful. Therefore, the only numerically invariant relation on a preference scale is the \textbf{\textit{k}-ratio}:
\begin{equation}
\label{eq: kfactor}
  k = \frac{p_a - p_b}{p_c - p_d}  
\end{equation}

where $p_a, p_b, p_c, p_d$ are preference scores forming two differences, whose ratio defines a relative scaling factor. Only affine transformations $p_{i} \mapsto a p_i + b$ (i.e., scaling and translation) applied uniformly to all preference scores preserve preference differences. Any other mathematical operation (e.g., absolute values or squaring) applied directly to individual scores is not meaningful under affine-scale invariance, because it alters difference ratios and thus changes the preference meaning. 

Linearity is a property of a vector space, not of human preferences themselves as such. Preference scores lie on a one-dimensional affine space, where only differences between scores carry meaning. Therefore, preference differences do not constitute \textbf{metric distances}, and preference aggregation \textbf{cannot} be based on metric distance measures or on operations applied directly to absolute numerical scores, a methodological issue widely noted in the MCDM literature \cite{french2023reflections}. Empirical studies confirm that many hybrid MCDM frameworks employing such operations exhibit structural inconsistencies, failing to ensure coherent or unique aggregated preference maximisation \cite{Ferdous2024}. This aligns with the broader recognition across the MCDM literature—spanning group decision‑making reviews \cite{Pajasmaa2025GroupDecision}, decision‑process reflections \cite{french2023reflections}, and methodological surveys \cite{figueira2016vol234}—that robust multi‑objective design and decision modelling requires tight integration between qualitative problem‑structuring and sound preference‑measurement theory. Classical representational measurement theory, as formulated by \cite{krantz1971foundations}, provides representational tools for physical quantities but does not identify the algebraic conditions under which human preference information can be subjected to meaningful mathematical operations. As Barzilai has shown \cite{barzilai2005}, this omission leads to the systematic misuse of addition, multiplication, and metric constructions in MCDM, because classical measurement theory does not establish the affine structure required for preference scales, where only differences—not absolute values—carry meaning. 

In contrast, Barzilai’s Preference Function Modelling (PFM) framework (\cite{barzilai2005}, 2006, 2010) explicitly provides these missing measurement‑theoretic foundations: it defines the formal conditions under which preferences can be represented and modelled through affine‑invariant operations, thereby resolving the deficiencies in the classical theory. For this PFM framework, Barzilai defined axioms\footnote{These are axioms in the formal mathematical sense—i.e., defining conditions that specify the admissible structure of preference representation.}, ensuring that preference scales and their combinations remain consistent and mathematically valid. These axioms specify the mathematical structure that preferences must satisfy to be meaningfully combined into a \textbf{single, consistent, and unique} (up to affine transformation) aggregated preference score, thereby establishing the formal conditions under which linear algebra and calculus (e.g. \cite{Strang2006}) can be rigorously applied to preference and other subjective quantities. For the purpose of preference aggregation in this work, the following four PFM-based axioms are outlined below. These form the theoretical foundation for pure MCDM.

A unique, consistent and well-defined MCDM outcome can be produced when preferences satisfy these axioms:

\vspace{-0.75em}

\subsubsection*{Axiom 1 – Preference Preservation ($\Delta$-meaningfulness)}\vspace{-0.5em}
Preference information is meaningful only in terms of \textbf{differences} between alternatives on an interval scale. 
These differences are relative, not absolute: distances as geometric constructs have an absolute zero, which does not exist in subjective preference measurement; zero-preference is always relative. 
Therefore, valid preference scales are only defined up to \textbf{affine transformations} (affine-invariant) applied uniformly:
\begin{equation}
\label{eq: affinetrans}
p_i' = a p_i + b, \qquad a > 0, \quad \forall i,
\end{equation}
which \textbf{preserve all meaningful preference information}. \textsc{Note:} Substituting Equation \eqref{eq: affinetrans} into Equation \eqref{eq: kfactor} leaves the $k$-ratio unchanged, thus preserving relative preference differences between alternatives.

\vspace{-0.75em}

\subsubsection*{Axiom 2 - Comparable Criteria}\vspace{-0.5em}
Preferences from different criteria can only be aggregated if they are measured on a common, valid, preference-based interval scale with commensurate units, ensuring that no criterion dominates due to its scale rather than its assigned weight. Thus, aggregation is only valid when \textbf{preference differences are commensurate}—i.e., measured in equivalent units across all criteria—ensuring that equal marginal differences carry equal preference meaning and no criterion dominance occurs.

\vspace{-0.75em}

\subsubsection*{Axiom 3 — Meaningful Zero-Reference}\vspace{-0.5em}
All criteria must share a \textbf{common, stable and meaningful zero-reference point (common affine reference origin)}, and all differences must be measured relative to this point. Only then, \textbf{linear aggregation} of these differences is mathematically valid and consistently interpretable across all alternatives and criteria; other non-linear aggregations — such as multiplicative, power-based, logarithmic, or distance-based optimisation methods — are therefore not meaningful under these conditions.

\vspace{-0.75em}

\subsubsection*{Axiom 4 — Uniqueness}\vspace{-0.5em}
Two preference systems producing identical judgments must correspond to the same underlying 
preference structure. Preference representation is therefore \textbf{unique up to affine transformations}, ensuring that criterion differences remain interval-invariant and equivalent preference information \textbf{cannot lead to conflicting aggregated rankings}.\\

Thus, whereas classical measurement theory does not establish decision-valid preference aggregation, PFM provides the mathematically rigorous measurement framework required for multi-criteria design–decision systems in MCDM. Yet, the operationalisation of aggregation within this framework remains to be established. In particular, an explicit aggregation operator is required that adheres to the four PFM axioms and yields a single, well-defined representative preference value per alternative. Only under such an operator does the alternative with the highest aggregated preference correspond to a unique best-fit-for-common-purpose solution, given the feasible set and the decision makers’ criteria and their relative importance. The development and justification of this PFM-consistent aggregation—governing mathematically valid and operational preference aggregation for MCDM—is the focus of this paper.

\section*{1. Constructing the Linear Preference Space (LPS)}

The objective of this section is to define a unified, interval-invariant scale of preference scores that fully complies with the PFM axioms, while supporting basic mathematical operations necessary for meaningful aggregation across criteria. For this, consider a set of alternatives $A_i$ ($i = 1,\dots,I$) evaluated against criteria $C_j$ ($j = 1,\dots,J$), with raw preference scores $p_{i,j}$ and associated non-negative weights $w_j$ ($w_j \ge 0$) such that $\sum_{j=1} w_j = 1$. 

To construct the \textbf{Linear Preference Space (LPS)}, the standard \textit{z}-score method is applied:
\begin{equation}
\label{eq: zscoring}
    z_{i,j} = \frac{p_{i,j}-\mu_j}{\sigma_j},
\end{equation}
where $\mu_j$ and $\sigma_j$ are the arithmetic mean and the standard deviation of the \textit{p}-scores for criterion $j$ respectively. Importantly, $\mu_j$ and $\sigma_j$ are purely arithmetic constructs used for centering and scaling; they carry no intrinsic preference meaning and so any mathematical operations may be applied to compute them. Only then the resulting $z_{i,j}$ values retain full preference meaning, maintaining all ratios of differences (\textit{k}-ratios) across alternatives and criteria. Explicitly the z-transformation is affine and can be expressed as  $z = a_j p + b_j$ with $a_j = 1/\sigma_j$ and $b_j = -\mu_j/\sigma_j$.

Moreover, by definition, \textit{z}-normalization ensures that, for each criterion $\mu_J=0$ and $\sigma_J=1$, where $\mu_J$ and $\sigma_J$ are the arithmetic mean and the standard deviation of the \textit{z}-scores respectively. These normalized values do not acquire intrinsic preference meaning; $\mu_J = 0$ serves as the arithmetic mean of the $z$-scores and provides a common and stable zero-reference point, so an alternative with $z_{i,j}=0$ can be interpreted as being ``average'' relative to the other alternatives on that criterion, while $\sigma_J = 1$ defines the unit of measurement (``standard deviation``) for the normalized differences, ensuring that preference differences are commensurate across all criteria.

In summary, the LPS construction is defined by the mapping
$T_j : p_{i,j} \mapsto z_{i,j}$, using $\mu_j$ and $\sigma_j$ to set a consistent origin and commensurate unit across all criteria. Z-normalization is a purely affine transformation of the raw scores that preserves all preference differences and establishes a common affine coordinate system. y fixing this common and stable origin across all criteria, the affine preference space obtains a vector-space representation, as every affine space becomes a vector space once an origin is chosen. This representation permits linear operations while preserving the affine structure and meaning of the underlying preference space. Moreover, this transformation provides a stable reference frame for aggregation fully consistent with the PFM axioms. The resulting $z$-scores retain the meaning of all relative preference differences while avoiding the notion of an absolute zero, ensuring that preferences remain represented as interval differences. After normalization, preferences are represented in a linear preference space where aggregation is meaningful only through affine-linear operations, since only such operators preserve the affine structure and the affine-invariant ratios between preference differences. Any non-linear operation would alter these ratios and thereby change the underlying preference meaning. Linearity is therefore a property of the constructed LPS, not of human preferences by themselves.

\subsection*{2. Aggregated Preference Ranking}

Given the same set of alternatives $A_i$, criteria $C_j$, criterion weights $w_j$, and normalized preference scores $z_{i,j}$ residing in the LPS as introduced in Section 1, the next step is to construct a single representative aggregated preference score $P_i^*$ for each alternative $A_i$. Formally, this representative score is defined as a unique scalar value that provides a best-fit to the weighted preference differences $z_{i,j}$ of alternative $A_i$ across the criteria. From a purely mathematical perspective, such a best-fit in a vector space is obtained by minimizing the total weighted least squared distance (WLSD) relative to a scalar representative value $F$, which yields the weighted centroid:

\begin{equation}
\label{eq:Fcentroid}
\min_{F} \sum_j w_j (f_{i,j} - F)^2
\quad \Longleftrightarrow \quad
F^* = \sum_{j} w_j \cdot f_{i,j}, \quad \text{where } \sum_{j} w_j = 1 .
\end{equation}

Here, the solution $F^*$ is the \textbf{linear weighted centroid} of the points $f_{i,j}$, obtained via a distance-based minimization. However, preferences are \textbf{differences} and \textbf{not distances}; therefore, preference aggregation cannot be derived from distance operations. In the context of preferences, where $F = P$ and $f_{i,j} = z_{i,j}$, the LPS is not a full vector space but an affine linear space of preference differences, a subset of a vector space (LPS $\subset V$), in which only linear operations on differences—scalar multiplication and addition—preserve preference meaning. Accordingly, while the WLSD formulation in Equation~\eqref{eq:Fcentroid} remains mathematically valid as an auxiliary derivation, neither the squared differences nor the (distance) minimization criterion as intermediate results carry preference meaning or ordering information. They serve solely as a mathematical device to expose the unique linear outcome admissible within the LPS. Hence, the only preference-theoretically valid result is its \textbf{linear solution}, the weighted centroid. 

\vspace{-0,8em}
\subsubsection*{Weighted centroid $P_i^*(z)$}\vspace{-0.5em}

Thus, in the LPS, the aggregated preference of an alternative is defined exclusively by linear operations on normalized preference differences. Accordingly, for an alternative $A_i$ with criterion scores $z_{i,j}$ and weights $w_j$, the representative aggregated preference score $\bm{P_i^*}$ is uniquely defined as the \textbf{weighted centroid} of its $z$-scores:

\begin{equation}
\label{eq: weightedcentroid}
P_i^*(z) = \sum_j w_j \cdot z_{i,j} \quad \text{where } \sum_{j} w_j = 1 
\end{equation}

This operator is a linear combination of preference differences (i.e., differences relative to the stable zero-reference point, obtained by centering the preference scores via z-normalization) and is therefore fully compatible with the affine preference structure of the LPS. The resulting $P_i^*(z)$ resides in the same affine space as the underlying $z$-scores, preserves all ratios of preference differences, and remains invariant under affine transformations of the original preference scales. \textsc{Note:} Within the multi-objective optimisation method IMAP (see see \cite{Wolfert2023}) this $P_i^*(z)$-based preference function aggregation operator is named \textit{A}, the \emph{a-fine-aggregator} (see \cite{Teuber2025Odycon}), which computes an aggregated preference score for candidate alternatives across multiple objectives.  

\vspace{-0,8em}
\subsubsection*{\textbf{$\boldsymbol{P}_i^*(z)$ Ranking}}\vspace{-0.5em}

Given the aggregated preference per alternative, the resulting \textbf{ranking} is obtained by ordering the alternatives according to their aggregated preference scores $P_i^*(z)$ in descending order, i.e.,
\[
\arg\max_i \, P_i^*(z).
\]
The alternative with the highest aggregated preference $P_i^*$ is therefore most preferred relative to the other alternatives and taking into account all criteria and their associated weights. For interpretability, it is often convenient to rescale the aggregated $P_i^*$ values to the interval $[0,100]$, so that the relatively worst (min: least preferred) and best (max: most preferred) alternatives correspond to 0 and 100, respectively. This can be achieved using the following linear min–max transformation, which preserves the relative ordering:

\begin{equation}
\label{eq:minmax}
P_i^{*\text{ scaled}} =
\frac{P_i^* - \min(P^*)}{\max(P^*) - \min(P^*)} \cdot 100,
\end{equation}

where $\min(P^*)$ and $\max(P^*)$ denote the minimum and maximum aggregated preference values across all alternatives. The endpoints 0 and 100 thus represent the relative minimum and maximum within the alternative set, with the zero point being purely relational rather than an absolute preference origin.

\textsc{Note:} The aggregated $P_i^*$ values may also be shifted or rescaled using any affine transformation (cf.\ Equation~\eqref{eq: affinetrans}) without affecting the ranking. For instance, the interval $[0,100]$ could equivalently be mapped to $[100,200]$, preserving both preference meaning and ordering. Because the min–max transformation is linear, it preserves all relative preference differences, keeping the \emph{k}-ratios from Equation~\eqref{eq: kfactor} fully consistent and thus maintaining affine invariance.

\newpage
\section*{3. Demonstrative Example }
Consider a $4 \times 3$ ranking example\footnote{In the Odesys book (see \cite{Wolfert2023}) similar MCDA examples can be found in Chapter~5, where the alternatives are referred to as variants. The example used here is taken from that Chapter.}
: i.e., four alternatives $A_i$ and three criteria $C_j$ with weights $w_j$, and 12 given preference scores $p_{ij}$. Here, $C_1$, $C_2$, and $C_3$ represent functionality, footprint, and costs, respectively. 

\[
\scalebox{0.95}{
$
\begin{array}{c|ccc}
\text{Alternative } A_i & C_1 (\text{functionality}) & C_2 (\text{footprint}) & C_3 (\text{costs}) \\ \hline
A_1 & 100 & 0   & 90  \\
A_2 & 0   & 100 & 100 \\
A_3 & 20  & 45  & 55  \\
A_4 & 85  & 60  & 0   \\ \hline
\text{Weights } w_j & 0.4 & 0.1 & 0.5\\
\end{array}
$
}
\]

The preference scores range from 0 to 100, where 0 represents the `worst' performance among the alternatives (i.e., the least preferred alternative relative to the others), and 100 represents the `best' performance among the alternatives (i.e., the most preferred alternative relative to the others). The value 0 therefore denotes the lowest relative score on the preference scale and does not imply the absence of performance on the corresponding criterion (i.e, no zero-performance). Our goal is now to compute a mathematically meaningful aggregated preference score for each alternative, allowing us to rank the alternatives from highest to lowest. To this end, the \(p\)-scores are first transformed into \(z\)-scores using Equation~\eqref{eq: zscoring}, yielding the corresponding \textbf{\(z\)-scores} for the four alternatives (\(A_1\)–\(A_4\)) and three criteria (\(C_1\)–\(C_3\)), as reported in the table below.

\[
\scalebox{0.95}{
$
\begin{array}{c|ccc}
\textit{z}\text{-scores} & C_1 & C_2 & C_3 \\ \hline
A_1 & 1.1556 & -1.4327 & 0.7351 \\
A_2 & -1.2148 & 1.3628 & 0.9908 \\
A_3 & -0.7408 & -0.1747 & -0.1598 \\
A_4 & 0.8000 & 0.2446 & -1.5660 \\ \hline
w_j & 0.4 & 0.1 & 0.5
\end{array}
$
}
\]

For these criteria, it can be verified : \textbf{(1)} that \(\mu_{C_1}=\mu_{C_2}=\mu_{C_3}=0\) and \(\sigma_{C_1}=\sigma_{C_2}=\sigma_{C_3}=1\); \textbf{(2)} using Equation~\eqref{eq: kfactor} that the \(k\)-factors in both \(p\)- and \(z\)-scores are affine invariant across the criteria: e.g., \(k_{C1}(p)=k_{C1}(z)=(A_{1}-A_{2})/(A_{4}-A_{3})\approx 1.53846\); \textbf{(3)} using Equation~\eqref{eq: affinetrans} that the \(p\)- and \(z\)-scores are affinely related and therefore affine invariant: e.g., for criterion \(C_1\) the \(p\)-values are mapped to the corresponding \(z\)-scores by \(z_{i,1} = 0.02370\,p_{i,1} - 1.21485\) (all numbers rounded to 5 decimal places for verification).

Now we can compute the aggregated $P_i^*$ per alternative using the weighted centroid from Equation~\eqref{eq: weightedcentroid} and the affinely invariant scaling to $[0,100]$ using Equation~\eqref{eq:minmax}. The final ranking is determined from highest to lowest $P_i^*$, resulting in the following table:

\begin{table}[h!]
\centering
\begin{tabular}{c c c c c}
\toprule
Alternative & $P_i^*$ & scaled [0-100] & Ranking \\
\midrule
$A_1$ & 0.68651 &  100 & 1 \\
$A_2$ & 0.14572 &  52 & 2 \\
$A_3$ & -0.39368 &  4 & 3 \\
$A_4$ & -0.43855 &  0 & 4 \\
\bottomrule
\end{tabular}
\label{tab:aggregated_ranking}
\end{table}

Lastly, the $z$-scores, along with their weight contributions and the aggregated preference $P_i^*$, are graphically depicted below in Figure 1. 

\begin{figure}[h!]
    \centering
    \includegraphics[width=0.85\linewidth]{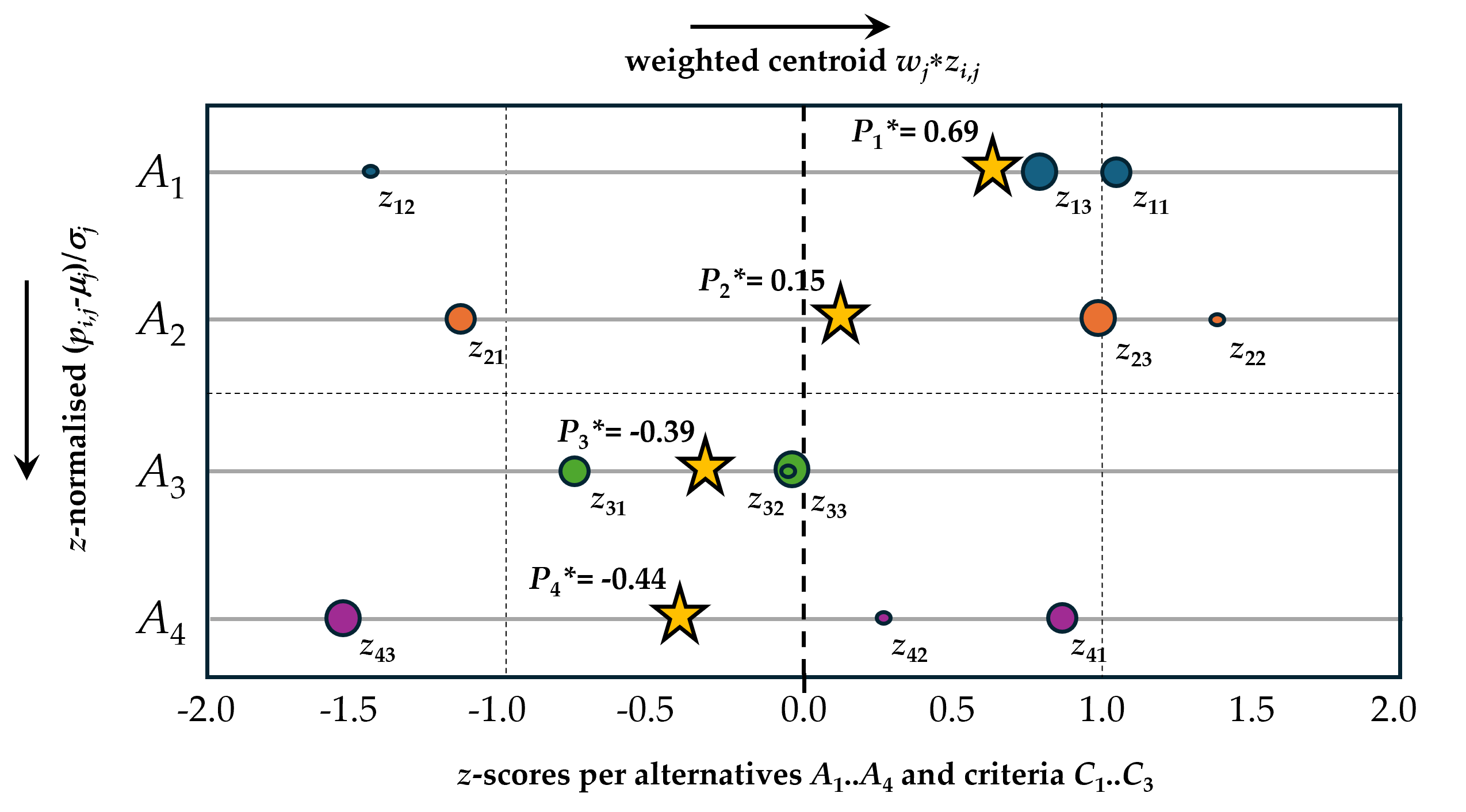}
    \vspace{-8pt}
    \caption{$P_i^*$ as aggregated preference per alternative}
    \label{fig:pstar_graph2}
\end{figure}

The figure illustrates the aggregated preference $P_i^*$ for each alternative, showing how the weighted criteria combine into a single representative value per alternative. Higher $P_i^*$ values indicate stronger aggregated preference and a better overall fit. This visual summary reflects the ranking in Table~\ref{tab:aggregated_ranking} and clearly illustrates the relative preference among alternatives, with Alternative 1 being the most preferred, as indicated by its position furthest to the right.

\subsubsection*{Intuitive Meaning}\vspace{-0.5em}

In the normalized LPS (defined by $z$-scores), $P_i^*$ can be interpreted as the \textbf{barycentre} of a set of points in an affine (linear) space. Graphically, this corresponds to the center of mass of the points, where each point contributes proportionally to its criterion weight. Intuitively, $P_i^*$ identifies the single point per alternative that best balances all weighted criterion preferences. \textsc{Note:} the center of mass per alternative across the weighted criteria, including $P_i^*$, equals zero (i.e., ``horizontal'' equilibrium), just as the centroid of the alternatives for each criterion also equals zero (i.e., ``vertical'' equilibrium).
Although preferences are \emph{not} physical distances and have no absolute zero, the barycentre (center of mass) analogy remains helpful. To illustrate this intuition, consider a horizontal beam with three point-masses placed at different positions, as shown in the Table below. 

\begin{table}[h!]
\centering
\begin{tabular}{c c c}
\hline
Mass $m_j$ & Weight & Position $x_j$ \\
\hline
1 & 0.40 & 1.1556 \\
2 & 0.10 & -1.4327 \\
3 & 0.50 & 0.7351 \\
\hline
\end{tabular}
\end{table}

Here, the position values $x_j$ correspond to the $z$-scores of Alternative 1 across the three criteria. The support must be located at the center of mass to keep the beam in equilibrium. Each point-mass represents a single preference-point contribution. The equilibrium position of the beam is given by $x_{\text{eq}} = (\sum_j m_j x_j) / (\sum_j m_j) \approx 0.6865$, which exactly equals the aggregated preference $P_1^*$. Even with a negative-position mass, the weighted positions combine into a single balance point. This illustrates the analogy between $P_1^*$ and $x_{\text{eq}}$: although not physical distances, these values provide an intuitive ``center of mass'' for multiple weighted preferences, just as $P_i^*$ represents the aggregated preference in the normalized linear preference space.

\textsc{Note:} The weighted center of mass is fully symmetric in the $z$-scores. Scores per alternative are horizontally centered, including the aggregated preference $P_i^*$, while scores per criterion are vertically centered. For example, for criterion \(C_1\) with four normalized scores and equal masses ($m_j = 1$), the centroid is $\bar{x} = \frac{1}{4} (1.1556 - 1.2148 - 0.7408 + 0.8000) = 0$, demonstrating that the weighted center of mass is symmetric: horizontally for alternative scores (including $P_i^*$) and vertically for criterion scores.

\subsubsection*{Weighted-mean and Distance-based rankings fail}\vspace{-0.5em}

Let us now repeat the 4×3 example from the start of this section, but using the following preferences scores $p_{ij}$. 

\[
\begin{array}{c|ccc}
\text{Alternative } A_i & C_1\ (\text{functionality}) & C_2\ (\text{footprint}) & C_3\ (\text{costs}) \\ \hline
A_1 & 80 & 0   & 24  \\
A_2 & 0   & 100 & 68 \\
A_3 & 30  & 45  & 55  \\
A_4 & 95  & 60  & 0   \\ \hline
\text{Weights } w_j & 0.4 & 0.1 & 0.5
\end{array}
\]

Following a similar approach as at the start of this section, the resulting $P_i^*(z)$ ranking is $A_2\,(100) \succ A_3\,(79) \succ A_1\,(32) \succ A_4\,(0)$. In addition, we consider the well-known absolute weighted arithmetic mean (WAM). For this example, WAM leads to a full tie, $A_1 = A_2 = A_3 = A_4\,(50)$, which is not congruent with the $P_i^*(z)$ ranking. Moreover, any affine rescaling of individual criteria (e.g., C1 or C3 on the $[0,100]$ scale) alters the WAM rankings, even changing the 'best' alternative, whereas the $P_i^*(z)$ ranking remains invariant. WAM aggregates alternatives by weighting absolute criterion scores without considering their relative positions within the set of alternatives. As a result, an extreme score on a single criterion can dominate the aggregation outcome, even if the alternative performs worse on multiple other, including heavily weighted, criteria. Like other distance-based rankings, WAM lacks a relative reference framework and therefore fails to produce unique, context-consistent, preference-stable rankings in the presence of conflicting criteria.

Using another weighted mean aggregation, the weighted geometric mean (WGM), the resulting ranking is $A_3\,(100) \succ A_1 = A_2 = A_4\,(0)$. Similar to WAM, WGM violates the PFM axioms, as its induced ranking depends on the choice of aggregation and the raw scores rather than on an affine-invariant preference structure. Consequently, neither WAM nor WGM provides a consistent or unique MCDM ranking in this case (see Appendix for further illustrative WAM/WGM examples).

It is beyond the scope of this section to provide the detailed calculations for the Euclidean and Manhattan distance-based rankings, as the approaches are outlined in the Appendix. Following this, the resulting distance-based rankings are: the Euclidean ranking ($D_i^E$) $A_3 \succ A_1 \succ A_2 \succ A_4$ and the Manhattan ranking ($D_i^M$) $A_3 \succ A_2 \succ A_1 \succ A_4$, demonstrating that non-linear distance-based aggregation fails to produce consistent rankings and also does not coincide with the unique $P_i^*(z)$ ranking. Thus, distance-based optimisations violate the PFM axioms, as their induced rankings depend on the choice of distance metric and non-linear transformations rather than being based on an affine-invariant preference structure. Consequently, distance-based methods cannot provide a unique or consistent MCDM ranking in this case (see Appendix for further illustrative examples). Last but not least, many MCDM techniques—including the distance‑based and weighted‑mean approaches illustrated here—are known to produce non‑unique, unstable, and method‑dependent results, a concern also emphasised in \cite{french2023reflections}.

\subsection*{4. Summary and Conclusions}

A unique and consistent decision outcome can be produced if and only if preferences adhere to the \textbf{PFM axioms}, which formalize preferences as measurable, subjective differences within a one-dimensional affine space (see \cite{barzilai2010} or 2022). These axioms ensure that preferences are represented as interval-invariant differences relative to a meaningful zero, comparable across criteria, and uniquely defined up to affine transformations. Consequently, aggregation preserves all meaningful relative differences, avoids scale- or criterion-induced distortions, and yields a consistent and unambiguous ranking of alternatives.

The aggregated preference $P^*(z)$, presented in this work, is \textbf{unique by construction}: it is the aggregation operator consistent with all four PFM axioms. Any aggregation that violates these axioms may fail to produce a coherent or meaningful preference-based ranking and is therefore not decision-valid. This uniqueness follows from the fact that the PFM axioms fully characterize the admissible aggregation operators, restricting them to affine linear functionals. Under normalization, this admissible class reduces to a single form: the weighted centroid of the z-normalized preference scores.

As preferences are \textbf{differences in a one-dimensional affine space}, not distance vectors in a vector space, aggregation cannot arise from distance-based operations. Distances—Euclidean, Manhattan, or others—are arbitrary and can produce infinitely many outcomes, failing to preserve preference meaning, ordering, and consistency (see the critique of mathematically invalid aggregations in \cite{french2023reflections}). For this reason, a Linear Preference Space (LPS) is required to be explicitly constructed. Within the LPS, defined by a stable zero-reference and commensurate units of preference differences across criteria, aggregation of z-normalized preferences via the linear weighted centroid preserves all ratios of differences, respects the affine structure, and remains invariant under affine transformations of the original scales. This ensures that the resulting ranking is \textbf{single, consistent, and unique}, providing a rigorous and preference-theoretically valid decision outcome. Note that the weighted least squares difference (WLSD) is used solely as a mathematical device to derive this linear ranking structure within the LPS; the squared differences and the minimization criterion themselves do not carry preference meaning. The resulting linear weighted centroid avoids these pitfalls and provides a fully consistent framework for preference-based decision-making.

In contrast, many commonly used MCDM methods—see \cite{kumar2021handbook} for an overview—do not guarantee adherence to the four PFM axioms. Whenever these axioms are not satisfied, the resulting multi-criteria rankings may become inconsistent, non-unique, and therefore decision-invalid. Moreover, many group decision-making approaches in multi-objective optimisation—particularly those based on Pareto-dominance sets, distance-to-ideal compromises, or interactive procedures—avoid explicit PFM-based preference aggregation altogether, as confirmed by recent systematic reviews such as \cite{french2023reflections} and \cite{Pajasmaa2025GroupDecision}. As a result, their outcomes may not be regarded as valid group decisions.

Even hybrid MCDM frameworks that attempt aggregation frequently exhibit structural inconsistencies, thereby failing to guarantee coherent and unique group-preference maximisation, as evidenced, for example, by the frameworks reviewed in \cite{Ferdous2024}, with related observations in \cite{Pajasmaa2025GroupDecision}. This explains why, despite decades of methodological development, such approaches cannot be considered fully consistent decision-making tools and instead function primarily as mathematical constructs, sometimes without explicitly incorporating preference information. While they may support exploration of the design and decision space and help identify non-dominated alternatives, they cannot guarantee the identification of a single, uniquely determined best-fit solution across all decision-makers’ objectives. Consequently, the potential for fully preference-consistent design and decision-making remains largely unexploited. The PFM-based preference aggregation approach presented in this paper enables the identification of a unique ranking and is explicitly operationalised for practical application in decision-making contexts.\\

\section*{APPENDIX }
This Appendix presents Examples \#E1–E3, which demonstrate specific conceptual pitfalls in preference aggregation, complementing the theoretical discussion in the main text. We first consider a decision analyst seeking to review PFM prior to commencing his new assignment.

\subsection*{Example \#E1 --- Why the Weighted Mean Fails in MCDM}

A decision analyst has to evaluate two job offers, \textbf{Job A} and \textbf{Job B}, based on two criteria: \textbf{growth opportunities (C1)} and \textbf{salary (C2)}. The criteria are weighted as $w_1 = 0.6$ for growth and $w_2 = 0.4$ for salary. The raw values for growth and salary per job are given in the table below.

\[
\begin{array}{c|cc}
 & C_1 \text{(growth)} & C_2 \text{(salary-€)} \\ \hline
A_1\ (\text{Job A}) & 15 & 50{,}000 \\
A_2\ (\text{Job B}) & 20 & 45{,}000 \\ \hline
\textit{w}_{1,2} & 0.6 & 0.4
\end{array}
\]

Using the values from the table above, the analyst computes the weighted average mean. With the original €-salaries, the result is \( A \succ B \). After discovering that the salaries were actually given in k\$, and converting them accordingly, the same weighted average mean now yields \( B \succ A \), even though the underlying decision problem has not changed.

Thus, the weighted average produces an infinite number of non-equivalent outcomes and 
therefore fails in decision-making. The analyst then decides to apply Preference Function 
Modelling (PFM) to determine whether he can make a unique and sound choice.

\subsubsection*{Unique Decision Making Using the Aggregated $P_i^*$} \vspace{-0.5em}

The decision analyst considers the following three approaches to the problem:

\begin{enumerate}
    \item Starting from the raw performance scores, either in € or in k\$, the values are directly normalized into $z$-scores, after which the aggregated preferences $P_i^*$ are computed.
       
    \item Starting by converting all criterion performance scores, whether expressed in € or in k\$ (which becomes irrelevant), into raw preference scores between 0 and 100 ('best'=100 and 'worst'=0), with $A_1$ (Job A) having scores $(0,\, 100)$ and $A_2$ (Job B) having scores $(100,\, 0)$.
\end{enumerate}

\textbf{Both approaches} produce the same normalized \textit{z}-score patterns and so the unique $P_i^*$:

\[
\begin{array}{c|cc|c}
 & C_1 & C_2 & P_i^* \\ \hline
A_1 & -1 & 1 & -0.2 \\
A_2 & 1 & -1 & 0.2 \\ \hline
\textit{w}_{1,2} & 0.6 & 0.4 & 
\end{array}
\]

Even though both approaches start from different numerical representations (€, k\$, or raw preference values), after \textit{z}-normalization they all lead to : $\boldsymbol{B \succ A}$. This demonstrates that the PFM / z-normalized centroid is \textbf{representation-invariant}, unlike the weighted average mean, which produced contradictory rankings. The decision analyst was very pleased, as at that moment the most important life decision for him could be made in a unique and transparent way.

\subsubsection*{Undefined Weighted Arithmetic or Geometric Mean} \vspace{-0.5em}

After making the proper job decision for Job B, he started working in the Data Science and Engineering (DSE) department. His first task was to analyse the following design evaluation problem, which involved two alternatives and two criteria, with the preference scores given as follows:

\[
\begin{array}{c|cc}
 & C_1 (\text{quality}) & C_2 (\text{safety}) \\ \hline
A_1 & 0.9 & 0.2 \\
A_2 & 0.7 & 0.6 \\ \hline
\textit{w}_{1,2} & 0.5 & 0.5
\end{array}
\]

The decision analyst calculated the $P_i^*$ using $z$-scores and concluded that both alternatives were equal (i.e., an indifferent choice: 50/50 and equal ranks). He went back to his manager and said he could not make a decision. The manager overruled him and instructed him to use the Weighted Arithmetic Mean (WAM) or the Weighted Geometric Mean (WGM)\footnote{ In general, WAM and WGM are widely used aggregation methods.The absolute weighted arithmetic mean (WAM) is defined as $\mathrm{WAM} = \sum_{j=1} w_j |p_j|$, and the weighted geometric mean (WGM) as $\mathrm{WGM} = \prod_{j=1} |p_j|^{\,w_j}$, where $p_j$ is the preference score on criterion $j$ and $w_j$ is the corresponding weight.}, as he should have been taught and gave him the following references: \cite{kumar2021handbook} and \cite{Krejci2018}. The decision analyst returned, disappointed in his first project, but did not give up and proceeded with the following analysis.
He used, as instructed, the Weighted Arithmetic Mean $\text{WAM}_i = w_1 \cdot x_{i1} + w_2 \cdot x_{i2}$ and the Weighted Geometric Mean $\text{WGM}_i = x_{i1}^{\,w_1} \cdot x_{i2}^{\,w_2}$, and also performed the $P_i^*(z)$ computation (see Equation \eqref{eq: weightedcentroid}). He constructed the following four Examples

For his storyline to his manager, he began his analysis with the following 2×2 problem — two alternatives and two criteria, C1 (quality) and C2 (safety) — as presented in the preferences table:

\[
\begin{array}{c|cc}
& C_1 & C_2 \\ \hline
A_1 & 0.90 & 0.35 \\
A_2 & 0.60 & 0.60 \\ \hline
w_j & w_1 & w_2
\end{array}
\]

For Example~1, the criteria weights are set as \((w_1, w_2) = (0.6, 0.4)\). He computed the following seemingly consistent ranking:

\begin{table}[h!]
\centering
\begin{tabular}{lcccc}
\hline
\textbf{$\boldsymbol{A}_{\mathbf{i}}$} & \textbf{WAM} & \textbf{WGM} & $\bm{P_i^*(z)}$ & \textbf{Rank} \\
\hline
A1 & 0.68 & 0.618 & 0.2 & 1 \\
A2 & 0.60 & 0.600 & -0.2 & 2 \\
\hline
\end{tabular}
\label{tab:summary}
\end{table}

However, he proceeded with the following. In Example~2/3, he considered two sets of weights for the criteria, \((w_1, w_2) = (0.54, 0.46)\) and \((w_1, w_2) = (0.48, 0.52)\), and used the same preference scores as in Example~1.

\begin{table}[h!]
\centering
\renewcommand{\arraystretch}{1} 
\begin{tabular}{c|c c c c}
\toprule
$\boldsymbol{w}_j$ & $\boldsymbol{A}_{\mathbf{i}}$ & \textbf{WAM} & \textbf{WGM} & $\bm{P_i^*(z)}$ \\
\midrule
\multirow{2}{*}{$\begin{matrix} w_1 = 0.54 \\ w_2 = 0.46 \end{matrix}$} 
                  & $\boldsymbol{A}_{1}$ & 1 & 2 & 1 \\
                  & $\boldsymbol{A}_{2}$ & 2 & 1 & 2 \\
\midrule
\multirow{2}{*}{$\begin{matrix} w_1 = 0.48 \\ w_2 = 0.52 \end{matrix}$} 
                  & $\boldsymbol{A}_{1}$ & 1 & 2 & 2 \\
                  & $\boldsymbol{A}_{2}$ & 2 & 1 & 1 \\
\bottomrule
\end{tabular}
\label{tab:ranks_only}
\end{table}

As he was able to demonstrate, both the WAM and WGM rankings failed, as they reversed the ranking order and were, at times, inconsistent with \(P_i^{*}(z)\). He then proceeded to analyse the following Example~4:

This again resulted in an ambiguous ranking for both WAM and WGM, and was not consistent with \(P_i^{*}(z)\):
\[
\begin{array}{c|cc}
& C_1 & C_2 \\\hline
A_1 & 0.90 & 0.20 \\
A_2 & 0.50 & 0.55 \\ \hline
\textit{w}_{1,2} & 0.5 & 0.5
\end{array}
\]

\begin{table}[h!]
\centering
\begin{tabular}{lccc}
\hline
\textbf{$\boldsymbol{A}_{\mathbf{i}}$
} & \textbf{WAM} & \textbf{WGM} & \textbf{$\bm{P_i^*(z)}$} \\
\hline
A1 & 1 & 2 & 1=2 \\
A2 & 2 & 1 & 1=2 \\
\hline
\end{tabular}
\label{tab:ranking_example2}
\end{table}

To completely impress his manager, he performed a final 4×3 Example 5, i.e., four alternatives and three criteria — C1 (technical quality), C2 (safety), and C3 (reliability) — with the following preference scores: 

\[
\begin{array}{c|ccc}
& C_1 (\text{quality}) & C_2 (\text{safety}) & C_3 (\text{reliability})\\ \hline
A_1 & 0.87 & 0.62 & 0.53 \\
A_2 & 0.80 & 0.72 & 0.72 \\
A_3 & 0.77 & 0.70 & 0.54 \\
A_4 & 0.79 & 0.93 & 0.56 \\ \hline
\textit{w}_{1,2,3} & 0.50 & 0.25 & 0.25
\end{array}
\]

A quick conclusion yields that WAM selects A4 as the best alternative, WGM selects A2 as the best alternative, and $P_i(z)$ selects A1 as the best alternative.

\begin{table}[h!]
\centering
\begin{tabular}{lccc}
\hline
\textbf{$\boldsymbol{A}_{\mathbf{i}}$
} & \textbf{WAM} & \textbf{WGM} & \textbf{$\bm{P_i^*(z)}$} \\
\hline
A1 & 4 & 3 & 1 \\
A2 & 2 & 1 & 2 \\
A3 & 3 & 4 & 4 \\
A4 & 1 & 2 & 3 \\
\hline
\end{tabular}
\label{tab:ranking_example4_final}
\end{table}

In this way, he demonstrated to his manager, using these examples, that the WGM and the WAM can yield infinitely many possible answers with different outcomes, non-consistent rankings, and are therefore both unsuitable as multi-criteria decision-making (MCDM) methods. Eventually, the roles were reversed: the manager was completely lost and had to abandon his incorrectly taught MCDM tools, while the decision analyst continued his work in the DSE department successfully. From that point on, they used only the PFM-based MCDM approach, as it was clear that all other methods (WAM and WGM) simply failed.

\subsection*{Example \#E2 --- Why we need Commensurate Criteria and a Stable Zero-Reference}
The decision analyst from the previous Example \#E2 now wants to evaluate five alternatives on two criteria ('5x2 problem'), using the following preference scores:

\[
\begin{array}{c|cc}
\text{\scriptsize dataset (0)} & C_1\,(\text{cost}) & C_2\,(\text{time}) \\ \hline
A_1 & 100 & 40 \\
A_2 & 0   & 60 \\
A_3 & 20  & 45 \\
A_4 & 85  & 50 \\
A_5 & 60  & 55 \\ \hline
\textit{w}_{1,2} & 0.5 & 0.5
\end{array}
\]

Having learned in Example \#E1 that directly applying \(P_i^{*}(z)\) yields an immediate, unique, and consistent ranking, he was---after the euphoria of that result---perhaps slightly overconfident. He began to wonder whether there might exist an alternative approach that would avoid \emph{z}-normalization altogether. In doing so, he decided first to explore the structure of the raw preference scales, hoping to deepen his understanding and, admittedly, to impress his manager even more than he already had in the previous Example.

Therefore, before applying any transformation, he examined whether the criteria scales were inherently comparable, noting that both \(C_1\) and \(C_2\) use \(0\text{--}100\) preference scores and share equal weights.

\subsubsection*{Scale Comparability}\vspace{-0.5em}
For a comparability check all criteria must at least have equal ranges and equal minimum values, i.e.\ \(\forall j,k \in \{1,\dots,m\}: p_{\min,j} = p_{\min,k} \ \wedge \ p_{\max,j} = p_{\max,k}\).

This condition ensures that a difference of \(\Delta p\) on any criterion contributes proportionally, independently of the original scale. This is exactly the essence of Barzilai's \emph{scale validity} (see \cite{barzilai2010}). By checking the above condition, it became clear that these scales are not directly comparable. Intuitively, he could also see that a difference of \(\Delta p = 10\) has very different implications for the two criteria: for C1 (cost), which ranges from 0 to 100, a difference of 10 represents 10\% of the total spread and is therefore substantial, while for C2 (time), which only varies from 40 to 60, a difference of 10 represents 50\% of the smaller spread; however, because the actual observed differences between alternatives on C2 are smaller (mostly 5--10 units), the criterion has less discriminative power in practice. Consequently, C2 contributes less to the overall ranking compared with C1.
Moreover, C1 represents costs ranging from €10,000 to €1,000,000 ('lower is better'), whereas C2 represents time varying only from 80 to 120 days ('shorter is better'). Therefore, C1 and C2 cannot be directly aggregated without proper normalization, because it is their spreads—not their weights—that determine their influence on the aggregation. This criterion dominance is not what weights represent; it is purely an artefact of the mapping, not true importance.  
Consequently, he decided to perform an affine transformation on C2, mapping all values to the interval \([0, 100]\).

\[
\begin{array}{c|cc}
\text{\scriptsize dataset (1)} & C_1 & C_2 \\ \hline
A_1 & 100 & 0   \\
A_2 & 0   & 100 \\
A_3 & 20  & 25  \\
A_4 & 85  & 50  \\
A_5 & 60  & 75  \\ \hline
w_{1,2} & 0.5 & 0.5
\end{array}
\]

Now it can be seen from Dataset (1) that the scales are comparable, since $p_{\min}=0$ and $p_{\max}=100$. After this analysis, the decision analyst then considered whether the $k$-ratio approach could still 'uniquely' determine a proper aggregated preference ranking. To address this question, he first compiled the following interlude.

\subsubsection*{Interlude}\vspace{-0.5em}
Following the basic \textit{k}-ratio (see Equation \eqref{eq: kfactor}), it can be shown that for multiple alternatives: let \(i = 1, \dots, n\) index the alternatives, and multiple criteria : let
\(j = 1, \dots, m\) index the criteria the \textit{k}-ratio of alternative \(i\) on criterion \(j\) is defined as

\begin{equation}
\label{eq:kratiogeneric}
k_{i,j} = \frac{p_{i,j} - p_{\min,j}}{p_{\max,j} - p_{\min,j}},
\end{equation}

where \(p_{\min,j}\) and \(p_{\max,j}\) denote the minimum and maximum preference scores for criterion \(j\). Then he assumed that a subset of the linear operations were valid and since the \textit{k}-scores are affine invariant, a consistent aggregation across criteria could be obtained using the following weighted centroid:

\begin{equation}
\label{eq: kratiocentroid}
\Pi_i^*(k) = \sum_{j=1} w_j \, k_{i,j}
\end{equation}

where $w_j$ is the weight of criterion $j$, with $\sum_{j=1} w_j = 1$. For $p_{\min,j}=0$ and $p_{\max,j}=100$, we obtain $\Pi_i^*(k)=\frac{1}{100} \sum_{j=1} w_j \, k_{i,j}$.

\subsubsection*{Comparison $\Pi_i^*(k)$ and $P_i^*(z)$ Ranking}\vspace{-0.5em}

Now he wanted to test and compare his new $k$-score--based $\Pi_i^*(k)$ ranking with the $P_i^*(z)$ ranking (see Equation \eqref{eq: weightedcentroid}). In addition to Dataset~(1), he therefore constructed the following additional Dataset (2), in which the preference values were slightly adjusted to become fully symmetric ([0,100]), and in which an extra intermediate reference was also introduced:

\[
\begin{array}{c|cc}
\text{\scriptsize dataset (2)} & C_1 & C_2 \\ \hline
A_1 & 0   & 100 \\
A_2 & 100 & 0   \\
A_3 & 50  & 50  \\
A_4 & 80  & 30  \\
A_5 & 20  & 70  \\ \hline
w_{1,2} & 0.5 & 0.5
\end{array}
\]

He then performed the computational analysis for both $P_i^*(z)$ and $\Pi_i^*(k)$, obtaining the following results:

\begin{table}[h!]
\centering
\begin{tabular}{c c!{\hskip 0.5em}|!{\hskip 0.5em}c}
\hline
\textbf{dataset} 
& $\boldsymbol{\Pi^*_i(k)}$ 
& $\boldsymbol{P^*_i(z)}$ \\
\hline
(1) 
& $A_4 = A_5 \succ A_1 = A_2 \succ A_3$ 
& $A_5 \succ A_4 \succ A_2 \succ A_1 \succ A_3$ \\
(2) 
& $A_4 \succ A_1 = A_2 = A_3 \succ A_5$ 
& $A_4 \succ A_1 \succ A_3 \succ A_2 \succ A_5$ \\
\hline
\end{tabular}
\end{table}

To his astonishment, he found that for both datasets he could not arrive at a unique ranking, as 
$\Pi_i^*(k)$ yielded different results than $P_i^*(z)$. He was particularly struck by the fact that 
$\Pi_i^*(k)$ contained several ties or indifferent outcomes, even for three alternatives, preventing 
a unique decision from being made. He then proceeded to re-examine the following 3x2-datasets to determine whether these ties were specific or indicative of a more general issue.

\vspace{-1.1em}
\begin{table}[h!]
\centering

\begin{minipage}{0.45\textwidth}
\centering
\[
\begin{array}{c|cc}
\text{\scriptsize dataset (10)} & C_1 & C_2 \\ \hline
A_1 & 0   & 0   \\
A_2 & 0   & 100 \\
A_3 & 100 & 60  \\ \hline
w_{1,2} & 0.30 & 0.70
\end{array}
\]
\end{minipage}
\hfill
\begin{minipage}{0.45\textwidth}
\centering
\[
\begin{array}{c|cc}
\text{\scriptsize dataset (15)} & C_1 & C_2 \\ \hline
A_1 & -50 & -50 \\
A_2 & -50 & 50  \\
A_3 & 50  & 10  \\ \hline
w_{1,2} & 0.30 & 0.70
\end{array}
\]
\end{minipage}

\vspace{0.5cm}

\begin{minipage}{0.45\textwidth}
\centering
\[
\begin{array}{c|cc}
\text{\scriptsize dataset (20)} & C_1 & C_2 \\ \hline
A_1 & 90  & 20  \\
A_2 & 50  & 55  \\
A_3 & 10  & 100 \\ \hline
w_{1,2} & 0.50 & 0.50
\end{array}
\]
\end{minipage}
\hfill
\begin{minipage}{0.45\textwidth}
\centering
\[
\begin{array}{c|cc}
\text{\scriptsize dataset (30)} & C_1 & C_2 \\ \hline
A_1 & 95  & 10  \\
A_2 & 60  & 55  \\
A_3 & 10  & 100 \\ \hline
w_{1,2} & 0.51 & 0.49
\end{array}
\]
\end{minipage}

\end{table}

After computing $P_i^*(z)$ and $\Pi_i^*(k)$ for the datasets (10), (15), (20), and (30), he obtained the following remarkable results:  

\begin{table}[h!]
\centering
\begin{tabular}{c c!{\hskip 0.5em}|!{\hskip 0.5em}c} 
\hline
\textbf{dataset} & $\boldsymbol{\Pi_i^*(k)}$ & $\boldsymbol{P_i^*(z)}$ \\
\hline
(10) & $A_3 \succ A_2 \succ A_1$ & $A_2 \succ A_3 \succ A_1$ \\
(15) & $A_3 \succ A_2 \succ A_1$ & $A_2 \succ A_3 \succ A_1$ \\
(20) & $A_1 = A_3 \succ A_2$     & $A_1 \succ A_3 \succ A_2$ \\
(30) & $A_1 \succ A_2 \succ A_3$ & $A_2 \succ A_1 \succ A_3$ \\
\hline
\end{tabular}
\end{table}

\vspace{-0.5em}
After careful analysis, the decision analyst concluded that $\Pi_i^*(k)$ could not yield a unique or consistent decision, whereas only $P_i^*(z)$ provided a stable ranking. Overall, $\Pi_i^*(k)$ fails to satisfy PFM Axioms because it lacks a common zero-reference and does not preserve linearity in meaningful differences, whereas $P_i^*(z)$ satisfies all axioms by design, yielding a unique and consistent ranking.

\textsc{Note:} \textit{k}-Ratio vs.\ \textit{z}-Normalisation — Preferences represent \textbf{differences}, not (vector) distances. They are relative, with no intrinsic absolute zero, and even \textit{z}-normalized preferences remain points in a linear preference space. The zero corresponds to the sample mean and a unit difference (standard deviation), serving only as a reference for linear operations. In a \textit{z}-normalized space, an aggregated ranking ${P_i^*(z)}$ can be computed with as few as two alternatives. In the \(p\)-domain (raw scores), at least three alternatives are generally needed to define meaningful \textit{k}-ratios (\(k_i = (p_i - p_{\min}) / (p_{\max} - p_{\min})\)), because no natural stable zero-reference exists. The \textit{z}-normalization introduces a functional zero at the mean (\(\mu_j\)) for each criterion, rendering all differences affine-invariant and allowing linear aggregation with only two alternatives. The functional mean can be interpreted as a “third” reference alternative in the \textit{k}-ratio framework.

\vspace{1em}
\subsection*{Example \#E3 --- Why Distance-based ranking fails in MCDM}

In this final example, the decision analyst aims to reveal a common pitfall in preference aggregation, namely the frequent—but incorrect—treatment of preferences as distances and the use of (non-linear) distance-optimisation--based rankings. To this end, he starts by evaluating a set of three alternatives and two criteria with the following \(p\)- or \(z\)-scores.

\[
\begin{array}{c|cc}
\text{\scriptsize \textit{p}-scores} & C_1 & C_2 \\ \hline
A & 90 & 10 \\
B & 50 & 50 \\
C & 10 & 90 \\ \hline
\textit{w}_{1,2} & 0.5 & 0.5
\end{array}
\qquad
\begin{array}{c|cc}
\text{\scriptsize \textit{z}-scores} & C_1 & C_2 \\ \hline
A & 1.0 & -1.0 \\
B & 0.0 & 0.0 \\
C & -1.0 & 1.0 \\ \hline
\textit{w}_{1,2} & 0.5 & 0.5
\end{array}
\]

\subsubsection*{Comparison of Distance-Based $D_i$ and $P_i^*$ Rankings}\vspace{-0.5em}

In any vector space, one may define a total weighted least squared difference (WLSD).
In the present context, this first takes the form of a weighted Euclidean distance (squared) between an alternative $i$ and an ideal (target) point $(1,1,\ldots,1)$, as commonly used in TOPSIS-style constructions (see, e.g., \cite{kumar2021handbook}). The normalized Euclidean
WLSD is then given by
\[
D_i^{E} = \sum_{j=1} w_j \,(z_{ij} - 1)^2,
\]
where $w_j$ are the criterion weights and $z_{ij}$ the normalized preference score of
alternative $i$ on criterion $j$. The resulting ranking is obtained by $\arg\min_{i} (D_i^{E})$, i.e., by minimizing the distance to the ideal point over all alternatives. The intention is to compare this ranking with the one induced by maximizing the aggregated preference score $P_i^*$ (see Equation \eqref{eq: weightedcentroid}). This resulted in:

\[
\begin{array}{c|cc|c|c}
 & C_1 & C_2 & \boldsymbol{P_i^*} & \boldsymbol{D_i^E} \\ \hline
A & 90 & 10 & 0.0 & 2.0 \\
B & 50 & 50 & 0.0 & 1.0 \\
C & 10 & 90 & 0.0 & 2.0 \\ \hline
\textit{w}_j & 0.5 & 0.5 & - & -
\end{array}
\]

He concluded, for the time being, that the Euclidean distance-based ranking 
\(D_i^{E}: B \succ A = C\) was inadequate, as it neither corresponded to his intuition—according to which a clear tie should emerge nor to the outcome of \( P_i^* \), which indeed indicated a full tie (as expected a-priori). As an initial reflection, however, he did not discard the distance-based ranking altogether, reasoning that the observed discrepancy might be attributable to the complete symmetry of the dataset. To gain further confidence, the analyst considered a second distance-based ranking to assess the sensitivity of the results to the choice of norm, namely the Manhattan distance. Analogous to the WLSD Euclidean-distance approach, the Manhatten-based ranking is obtained by:

\[
\arg\min_{i}\, (D_i^{M}),
\quad \text{with} \quad
D_i^{M} = \sum_{j=1} w_j \, \lvert z_{ij} - 1 \rvert .
\]

In this second step, he applied the Manhattan-based WLSD to a set of non-symmetrical datasets for further evaluation, together with their corresponding \( p \)-scores.

\[
\begin{array}{c|cc}
\text{\scriptsize dataset (1)} & C_1 & C_2 \\ \hline
A & 60 & 0 \\
B & 0 & 100 \\
C & 100 & 20 \\ \hline
\textit{w}_{1,2} & 0.3 & 0.7
\end{array}
\qquad
\begin{array}{c|cc}
\text{\scriptsize dataset (2)} & C_1 & C_2 \\ \hline
A & 90 & 20 \\
B & 50 & 55 \\
C & 10 & 100 \\ \hline
\textit{w}_{1,2} & 0.5 & 0.5
\end{array}
\qquad
\begin{array}{c|cc}
\text{\scriptsize dataset (3)} & C_1 & C_2 \\ \hline
A & 0 & 0 \\
B & 0 & 100 \\
C & 100 & 60 \\ \hline
\textit{w}_{1,2} & 0.3 & 0.7
\end{array}
\]

To his astonishment, he obtained the following results: all distance-based rankings—whether Euclidean or Manhattan (with ``lowest is best'')—differed from the aggregated $P_i^*$ preference ranking (with ``highest is best''). Moreover, the distance-based rankings could even differ among themselves. \textsc{Note:} all values reported in the tables are rounded. In essence, this demonstrates that distance-based optimisation can lead to inconsistent rankings and therefore cannot guarantee a unique ranking in MCDM. Such approaches fail to comply with the axioms of Preference Function Modelling (PFM), which are required for consistent and unique preference aggregation.

\begin{table}[h!]
\centering
\renewcommand{\arraystretch}{1.2} 

\begin{minipage}{0.45\textwidth}
\centering
\begin{tabular}{c|c|c}
\hline
\text{\scriptsize dataset (1)} & $\boldsymbol{P_i^*}$ & $\boldsymbol{D_i^E}$ \\
\hline
A & -0.599  & 2.805 \\
B & 0.582   & 1.690 \\
C & 0.017   & 1.505 \\
\hline
\scriptsize ranking & $B \succ C \succ A$ & $C \succ B \succ A$ \\
\hline
\end{tabular}
\end{minipage}
\hfill
\begin{minipage}{0.45\textwidth}
\centering
\begin{tabular}{c|c|c}
\hline
\text{\scriptsize dataset (2)} & $\boldsymbol{P_i^*}$ & $\boldsymbol{D_i^M}$ \\
\hline
A & 0.027  & 1.197 \\
B & -0.051  & 1.051 \\
C & 0.024  & 1.248 \\
\hline
\scriptsize ranking & $A \succ C \succ B$ & $B \succ A \succ C$ \\
\hline
\end{tabular}
\end{minipage}

\vspace{0.5cm} 

\begin{minipage}{0.6\textwidth}
\centering
\begin{tabular}{c|c|c}
\hline
\text{\scriptsize dataset (3)} & $\boldsymbol{D_i^E}$ & $\boldsymbol{D_i^M}$ \\
\hline
A & 4.571  & 2.121 \\
B & 0.887 & 0.607 \\
C & 0.542 & 0.710 \\
\hline
\scriptsize ranking & $C \succ B \succ A$ & $B \succ C \succ A$ \\
\hline
\end{tabular}
\end{minipage}

\end{table}

After the previous evaluations, he realised that when it comes to a life-or-death decision 
(e.g., choosing an operation for person A, B, or C), it is of paramount importance that the ranking is unique and consistent. 
The result, however, was non-unique under both the Euclidean and Manhattan distances, 
demonstrating that distance-based rankings—whether Euclidean, Manhattan, or based on other norms—are arbitrary 
and cannot provide a proper ranking of alternatives. 
Finally, he concluded that even within distance-based rankings, an infinite number of different outcomes can occur, 
and that \textbf{preferences are not (metric) distances}, but points in a one-dimensional affine space, with relative meaning.
Overall, it is concluded that distance-based rankings are norm-dependent and unstable, whereas the weighted centroid $P_i^*(z)$ in a linear preference space is unique by design, as confirmed by step-by-step numerical examples. This can be attributed to at least the following:

\begin{enumerate}\setlength{\itemsep}{0pt}\setlength{\parskip}{0pt}\setlength{\parsep}{0pt}
    \item $P_i^*$ is additive, linear, and therefore \textbf{consistent} with the preference axioms governing a linear preference space (LPS).
    \item Distance-based measures $D_i$ construct rankings using absolute (Manhattan) or squared (Euclidean) deviations, which are strictly non-negative by definition. As a result, they cannot be congruent with the signed centroid aggregation $P_i^*$, representing a linear preference projection, and may \textbf{directly contradict the preference ranking}. Any coincidence in ranking is incidental rather than structural.
    \item Distance-based aggregation relies on \textbf{an absolute metric zero}, which does not exist in preference scales, rather than on a meaningful preference centre. In contrast, $P_i^*$ is defined within a centred linear preference space, where the mean provides a functional zero and preserves all relative differences.
\end{enumerate}


\begin{thebibliography}{99}

\bibitem{Pajasmaa2025GroupDecision}
Pajasmaa, J., Miettinen, K., \& Silvennoinen, J. (2025).
Group decision making in multiobjective optimization: A systematic literature review.
\textit{Group Decision and Negotiation}, 34(2), 329--371.
https://doi.org/10.1007/s10726-024-09915-8

\bibitem{teuber2025odycon}
Teuber, L.G., van Heukelum, H.J., \& Wolfert, A.R.M. (2025).
Advancing strategic planning and dynamic control of complex projects.
\textit{arXiv}, https://doi.org/10.48550/arXiv.2408.12422.

\bibitem{Ferdous2024}
Ferdous, J., Bensebaa, F., Milani, A. S., Hewage, K., Bhowmik, P., \& Pelletier, N. (2024). 
Development of a generic decision tree for the integration of multi-criteria decision-making (MCDM) and multi-objective optimization (MOO) methods under uncertainty to facilitate sustainability assessment: A methodical review. 
\textit{Sustainability, 16}(7), 2684. https://doi.org/10.3390/su16072684

\bibitem{french2023reflections}
French, S. (2023). Reflections on 50 years of MCDM: Issues and future research needs. \textit{EURO Journal on Decision Processes}, 11-100030, https://doi.org/10.1016/j.ejdp.2023.100030

\bibitem{wolfert2023}
Wolfert, A.R.M. (2023).  \textit{Open Design Systems}. IOS Press/ TUDelft-OPEN.
https://doi.org/10.3233/RIDS10 or https://www.odesys.nl or
\\ https://books.open.tudelft.nl/home/catalog/book/78. 

\bibitem{barzilai2022}
Barzilai, J. (2022). \textit{Pure Economics}. FriesenPress.

\bibitem{kumar2021handbook}
Kumar, U.,  Chan, F.T.S. (Eds.). (2021). \textit{A Handbook on Multi-Attribute Decision-Making Methods}. Wiley.

\bibitem{krejci2018}
Krejčí, J., \& Stoklasa, J. (2018). Aggregation in the analytic hierarchy process: Why weighted geometric mean should be used instead of weighted arithmetic mean. \textit{Expert Systems with Applications, 114}, 97--106. https://doi.org/10.1016/j.eswa.2018.06.060

\bibitem{figueira2016vol234}
Figueira, J. R., Greco, S.,  Ehrgott, M. (Eds.). (2016). \textit{Multiple Criteria Decision Analysis: State of the Art Surveys} (Vol. 233 \& 234). Springer.

\bibitem{barzilai2010}
Barzilai, J. (2010). Preference function modelling: The mathematical foundations of decision theory. In M. Ehrgott, J. R. Figueira, \& S. Greco (Eds.), \textit{Trends in Multiple Criteria Decision Analysis} (pp. 57--86). Springer. https://doi.org/10.1007/978-1-4419-5904-1\_3

\bibitem{barzilai2006}
Barzilai, J. (2006). Preference modeling in engineering design. In K. E. Lewis, W. Chen, \& L. C. Schmidt (Eds.), \textit{Decision Making in Engineering Design} (pp. 43--47). ASME Press.

\bibitem{strang2006}
Strang, G. (2006). \textit{Linear algebra and its applications} (4th ed.). Cengage Learning.

\bibitem{barzilai2005}
Barzilai, J. (2005). Measurement and preference function modelling. \textit{International Transactions in Operational Research, 12}(2), 173--183. https://doi.org/10.1111/j.1475-3995.2005.00496.x

\bibitem{krantz1971foundations}
Krantz, D. H., Luce, R. D., Suppes, P., and Tversky, A. (1990).
\textit{Foundations of Measurement: Volumes I--III}.
Academic Press, New York.

\end{thebibliography}
\end{document}